\documentclass[10pt]{article}  
\usepackage{amsmath}
\usepackage{amssymb}
\usepackage{epsfig}
\usepackage{graphics}
\usepackage{theorem}

\theoremheaderfont{\scshape}

\theorembodyfont{\normalfont}

\numberwithin{equation}{section}
\newcommand{\bint}{\int \hspace{-11pt}{\bf\diagup} \hspace{4pt}}
\newcommand{\bintsml}
{\int \hspace{-7.5pt}{{\scriptscriptstyle\bf\diagup}} \hspace{4pt}}

\newcommand{\nit}{\mathbb{N}}
\newcommand{\rit}{\mathbb{R}}

\newcommand{\eps}{{\varepsilon}}
\newcommand{\avec}{\mathbf{a}}
\newcommand{\aavec}{\mathbf{A}}
\newcommand{\bvec}{\mathbf{b}}
\newcommand{\bbvec}{\mathbf{B}}
\newcommand{\tc}{{t}}

\newcommand{\xvec}{{\mathbf{x}}} 
\newcommand{\xxvec}{\mathbf{X}}
\newcommand{\yvec}{\mathbf{y}}
\newcommand{\yze}{\yvec^0}

\newcommand{\yyvec}{\mathbf{Y}}

\newcommand{\uuvec}{\mathbf{U}}
\newcommand{\vvec}{{\mathbf{v}}}
\newcommand{\vvvec}{\mathbf{V}}
\newcommand{\wvec}{\mathbf{w}}

\newcommand{\vtot}{\mathbf{m}}
\newcommand{\vmr}{ \mathbf{M}}

\newcommand{\pvmr}{\mathbf{N}}

\newcommand{\vwn}{\mathbf{W}}

\newcommand{\tov}{{\overline{t}}}
\newcommand{\ttov}{{\overline{T}}}

\newcommand{\lov}{{\overline{l}}}
\newcommand{\vov}{{\overline{v}}}

\newcommand{\ds}{\displaystyle}
\newcommand{\toe}{\frac{t}{\eps}}
\newcommand{\tmsoe}{\frac{t}{\eps}}

\newcommand{\tsep}{\frac{t}{\eps}}
\newcommand{\unsep}{\frac{1}{\eps}}

\newcommand{\fracp}[2]{\frac{\partial #1}{\partial #2}}

\newcommand{\xxze}{{\xxvec^0}}

\newcommand{\xxun}{{\xxvec^1}}

\newcommand{\yyze}{{\yyvec^0}}

\newcommand{\yyun}{{\yyvec^1}}

\title{Long term object drift forecast in the ocean with tide and wind}
\author{P. Ailliot$^\textbf{b}$ 
\and E. Fr\'enod$^\textbf{a,c}$ 
\and  V. Monbet$^\textbf{a,b,d}$}
\begin{document}
\maketitle

\begin{center}
$^\textbf{a}$ Laboratoire d'\'Etude et Mod\'elisation des Environnements Littoraux (Lemel), 
Universit\'e de Bretagne Sud, Centre Yves Coppens, Campus de Tohannic,
F-56000, Vannes.

$^\textbf{b}$ Laboratoire de Statistiques Appliqu\'ees de BREtagne Sud (Sabres), 
Universit\'e de Bretagne Sud, Centre Yves Coppens, Campus de Tohannic,
F-56000, Vannes.

$^\textbf{c}$ Laboratoire de Math\'ematiques et Applications des Math\'ematiques (LMAM) ,
 Universit\'e de Bretagne Sud, Centre Yves Coppens, Campus de Tohannic,
F-56000, Vannes.

$^\textbf{d}$ 
Hydrodynamique et Oc\'eano-m\'et\'eo (DCB/ERT/HO),
Ifremer, Centre de Brest, F-29280  Plouzan\'e.
\end{center}
                      

{\bf Abstract}
 In this paper, we propose a new method to forecast the drift of objects in 
near coastal ocean on a period of several weeks. The proposed approach
consists in estimating the probability of events linked to the drift
using Monte Carlo simulations.
It couples an averaging method which permits to decrease the 
computational cost and a statistical method in order to take into 
account the variability of meteorological loading factors.

~
{\bf Keywords} Averaging method, Two scale modeling, Monte Carlo simulation, Wind time series, Ocean forecasting, Object drift

~
\section{Introduction}

Drift of things in the ocean is potentially dangerous for human 
activities and marine ecosystems. For instance, drifting containers 
may cause serious accidents in the event of collision with ships 
and oil spills may have very negative impacts especially in coastal areas. 
In this paper, we investigate a new method to forecast the drift 
of an object in 
near coastal ocean on a period of several weeks. 

~

The motion of a drifting object on the sea surface is the net result 
of a number of forces acting upon it (water currents due to 
tide wave, atmospheric wind, wave motion,  wave induced currents, 
gravitational force and buoyancy force). It is possible 
to estimate the drift trajectory given information on the local 
wind, the surface current, and the shape and the buoyancy of the object.
For instance, in order to estimate the position of lost containers, the 
safety-and-rescue services generally use short-term meteorological 
forecasts as forcing of an hydrodynamic model of drift 
(see Daniel {\it et al.} \cite{DaJa2002}). It is usual in such 
problems to consider several possible buoyancy and drift properties 
for the object since these features are not known precisely in 
most cases. An uncertainty on the initial conditions (position and time) 
may also be taken into account (see Hackett {\it et al.} \cite{HaBr2004}). 

Since we are interested in longer periods of time in the present study, 
we cannot directly use meteorological forecasts to estimate the 
object trajectory. Then the drift forecast has to be led out
in terms of probability. This will be done using a Monte Carlo method, 
which makes it possible to estimate the probability of some scenarios linked
to the object's trajectory, like the probability of being in a given 
point at a given time or the
probability of running aground in given areas, for example. It 
consists in computing the object trajectories corresponding to a large number
of meteorological time series representative of the climatology on the 
considered area.
Because of the 
variability of the meteorological conditions, it is necessary to 
compute a large number of trajectories in order to get reliable 
estimates of the quantities of interest. Hence we are faced with two problems.

At first, as the existing data sets describe the meteorological 
conditions only on the few last decades, it is necessary to be 
able to simulate new realistic meteorological time-series. 
For this, we have used stochastic models. More precisely, 
it is assumed that these time series can be decomposed as 
the sum of two components, the first one that describes the 
meteorological conditions at a synoptic scale and the second 
one that represents fluctuations at the smaller scale. 
Then, the synoptic component is simulated with a non-parametric 
resampling algorithm proposed by Monbet, Ailliot and Prevosto 
\cite{MonbetAlliotPrevosto}, and the short-term component with 
an Autoregressive model. 

Secondly, we have to compute the corresponding 
object trajectories, which oscillate
with tide, in a reasonable computational time. 
Indeed, the trajectory of a drifting object submitted to tide wave currents and 
currents generated by other meteorological factors, is essentially 
composed by a trend and oscillations due to the tide,
which have a small period with respect to the time period of interest 
(several weeks).  Due to these high frequency oscillations, the 
numerical integration of the considered system is generally time consuming. 
So that it is convenient to split the estimation of the displacement 
of the object in two steps: first, we calculate the trend 
which is of low computational cost thanks to the lack of oscillations 
and then we reconstruct the oscillation around the trend. For this, 
we apply the Averaging Method developed by Fr\'enod \cite{frenod2004} 
in order to 
identify the averaged fields governing the trend and the 
oscillating operators allowing the reconstruction of the real object 
trajectory from the trend. As far as we know, this method has never been used 
before in the metocean field. 

~

In order to check the validity of the proposed methodology, 
we consider a situation with deterministic currents and negligible wave effects 
and we study a simplified model where the acceleration of the floating 
object is equal to the sum of the acceleration of the water (due to 
the wave tide and a perturbation of smaller order which represents 
other water currents) and the difference between the wind speed and 
the object velocity. 
This model is described more precisely in the second section. Then, 
in the third section, we present
the basic ideas of the Averaging Method and we compare computational cost of this method with the one which consists 
in integrating the real trajectory directly.
In Section 4, we validate our method using numerical experiments. 
First, the models of the meteorological forcing fields 
(tide wave, perturbation and wind) are specified. Then, the accuracy of 
the Averaging Method is checked on the basis of numerical comparisons. 
Finally, in Section 5 the methodology is illustrated with an example in which
the probability of running aground is estimated for an academic domain.

\section{Model} 
The model on which we shall implement the method evoked above 
is a very simplified model for large time drift
in ocean, above continental shelf in strong tide zone, of an almost
completely submerged object submitted to wind.
This model is extracted from exhaustive scaling analysis
for large time floating object drift which we shall present in 
a forthcoming paper.

~

In short, the evolution of the position  $\xxvec(\tc)=\xxvec(\tc;\xvec,\vvec)$ $\in\rit^2$
and the velocity $\vvvec(\tc)=\vvvec(\tc;\xvec,\vvec)$  $\in\rit^2$ of
the considered object, having $\xvec$ and $\vvec$ as initial position
and velocity, is given by:
\begin{align} \ds
\frac{d\xxvec}{d\tc} &=\vvvec, 
\label{EqPosOb}
\\  \ds
\frac{d\vvvec}{d\tc} &=\frac{d}{dt}\Big[\vtot(\tc,\xxvec)) \Big]
+\Lambda (\wvec(\tc,\xxvec) -\vvvec )
 =\fracp{\vtot}{\tc}(\tc,\xxvec)+ 
\big(\nabla\vtot(\tc,\xxvec)\big)\vvvec+ \Lambda (\wvec(\tc,\xxvec) -\vvvec ),
\label{EqVitOb}
\end{align}
where 
$\vtot \equiv\vtot(\tc,\xvec)$ is the ocean velocity field and 
$\wvec\equiv\wvec(\tc,\xvec)$ the wind velocity field.
$\nabla \vtot$ stands for the Jacobian matrix of $\vtot$.
This equation says nothing but that the object is submitted
to the sea water acceleration and to the wind force quantified by a constant $\Lambda$.

The time scale $\tov$ on which we want to observe the drift phenomenon 
is about 3 months and the object is submitted to tide oscillations
whose period $\ttov$ is about 12.5 hours. Then a small parameter $\eps$
appears in our problem: the ratio $\frac{\ttov}{\tov}$ tide period on 
observation time scale.
The magnitude of $\eps$ is about $1/200$.
We consider that the velocity measurement of the
object, wind and water are all done with the same unit: $\vov$.
The observation length scale $\lov$, which has to be the characteristic 
length of continental shelf, is about several hundred kilometers. 
This length has to be compared with the characteristic distance 
$\ttov \vov$ the water covers during a tide period, which is of
some kilometers. It seems then reasonable to consider that  
the ratio $\frac{\ttov \vov}{\lov}$ is also about $\eps$.

Having those scale considerations at hand, we introduce the following
rescaled variables  $t'$, $\xvec'$ and  $\vvec'$ expressing time, position 
and velocity in unit $\tov$, $\lov$ and  $\vov$
\begin{gather}
t=\tov t', \xvec =\lov\xvec'
\text{ and } \vvec=\vov\vvec',
\end{gather}
the rescaled trajectory $(\xxvec'(t';\xvec',\vvec')$, 
$\vvvec'(t';\xvec',\vvec'))$ defined by
\begin{gather}
\lov \xxvec'(t';\xvec',\vvec')= 
\xxvec(\tov t';\lov\xvec',\vov\vvec'), ~~~
\vov \vvvec'(t';\xvec',\vvec')= 
\vvvec(\tov t';\lov\xvec',\vov\vvec'),
\end{gather}
and the rescaled fields $\vtot'$ and $\wvec'$ defined by
\begin{gather}
\vov\vtot'(t',\xvec') = \vtot(\tov t',\lov\xvec'), ~~~
\vov\wvec'(t',\xvec') =  \wvec(\tov t',\lov\xvec').
\end{gather}

We consider that the sea velocity writes
\begin{gather}
\vtot'(t',\xvec') 
= \vmr(t',\frac{\tov}{\ttov} t',\xvec')
+ \eps \pvmr(t',\frac{\tov}{\ttov} t',\xvec')
= \vmr(t',\frac{t'}{\eps},\xvec')
+ \eps \pvmr(t',\frac{t'}{\eps},\xvec'),
\end{gather}
where 
$\vmr (\tc',\frac{t'}{\eps},\xvec')$  $\in\rit^2$ is the rescaled sea 
water velocity exclusively due to the tide wave.
It is supposed to be regular. 
Considering its dependency  with respect to the oscillating
time variable, we suppose that 
$\theta'\mapsto\vmr(\tc',\theta',\xvec')$
is a $1-$periodic function satisfying 
$\ds \bint \vmr (\tc',\theta',\xvec') \, d\theta' =0$ where
$\ds \bint \vmr (\tc',\theta',\xvec') \, d\theta' =
\int_0^1\vmr (\tc',\theta',\xvec') \, d\theta' $.
The field $\eps \pvmr (\tc',\frac{t'}{\eps},\xvec')$, where
$\pvmr (\tc',\theta',\xvec')$ is also $1-$periodic in $\theta'$,
is the sea water velocity perturbation induced by  
meteorological factors.
As concerns the wind velocity field, we consider that $\wvec'$
also involves two time scales, i.e.,
\begin{gather}
\wvec'(\tc',\xvec') = \vwn(\tc',\frac{t'}{\eps},\xvec').
\end{gather}
Nonetheless, wind time series
bring out as unrealistic considering a periodic dependency
of $\vwn(\tc',\theta',\xvec')$ with respect to $\theta'$.
In practice, we only consider that, for any $\tc'$  and $\xvec'$,
it admits an average value $\ds \bint \vwn (\tc',\theta',\xvec') \, d\theta'$ which actual definition is discussed later on.

Finally, we may deduce the equation satisfied by the rescaled trajectory
and involving the rescaled fields. Removing the $'$, 
this equation, which is the model on which we shall implement our
method, reads: 
\begin{align} \ds
&\frac{d\xxvec}{d\tc} =\vvvec, 
\label{SM1ts1}
\\  \ds
&
\begin{aligned}
  \frac{d\vvvec}{d\tc} =\unsep \fracp{\vmr}{\theta}(\tc,\tsep,\xxvec) 
  + \fracp{\vmr}{\tc}(\tc,\tsep,\xxvec)
  + \fracp{\pvmr}{\theta}(\tc,\tsep,\xxvec)
  + \big( \nabla \vmr (\tc,\tsep,\xxvec)\big)\vvvec 
  + \vwn(\tc,\tsep,\xxvec) -\vvvec ~~~~&
\\ 
  +\eps \Big(\fracp{\pvmr}{\tc}(\tc,\tsep,\xxvec)
  +\big( \nabla \pvmr (\tc,\tsep,\xxvec)\big)\vvvec\Big).&
\end{aligned}
\label{SM1ts2}
\end{align} 
Notice that system (\ref{SM1ts1})-(\ref{SM1ts2}) is rescaled, and that, in it, 
every variable and field characteristic scale is of order 1.

~

Clearly, the model under consideration is too simplistic to be
used for operational applications.
Nevertheless, it contains most of the physical ingredients
of drift of object in the ocean: joint action of  sea and 
wind, two time scales, possibility of using not so unrealistic
sea velocity fields.
Moreover, it seems to be relatively straightforward to incorporate a realistic sea velocity
in it, with a tide period that weakly evolves
with time and to apply it in a real geographical geometry.
Hence, the validity of the methodology we present in this paper
is not limited by the simplifications we consider.

\section{Asymptotic analysis}
Having the goal of using the Monte Carlo Method 
to estimate the probability of events linked to the trajectory of 
the drifting object in mind, we need
to compute the object trajectory for a large number of wind conditions.
The solution $(\xxvec,\vvvec)$ contains $\frac{1}{\eps}-$frequency
oscillations. Then, solving (\ref{SM1ts1})-(\ref{SM1ts2}) directly by 
numerical methods forces the use of a very small time step.
For instance, if the explicit Euler scheme is used, since the characteristic
size of the
left hand side of (\ref{SM1ts2}) and of its gradient are
about $\frac{1}{\eps}$ and the size of its time derivative is about 
$\frac{1}{\eps^2}$,
the classical error estimate yields, for a time step
$\Delta t$ small enough, an error about
\begin{equation}
\frac{\Delta t}{\eps} \bigg(1+\frac{\Delta t}{\eps}\bigg) ^{\frac{1}{\Delta t}}
\sim ~ \frac{\Delta t }{\eps}e^\frac{1}{\eps}.
\end{equation}
Hence, if we want to obtain a precision about $\eps^2$ a time step $\Delta t$
about $\eps^3e^{-\frac{1}{\eps}}$ is needed.
This is really too small for operational applications.
Hence we shall write an expansion of $(\xxvec,\vvvec)$ and find
non oscillating equations satisfied by the terms of this expansion.
This way, the previously evoked constraint imposed on the time 
step vanishes. 

~

It is an easy game to see that system (\ref{SM1ts1})-(\ref{SM1ts2})
enters the framework of an oscillatory-singularly perturbed 
dynamical system
\begin{equation}\label{Rsd}
\frac{d}{dt}\begin{pmatrix}\xxvec \\ \vvvec \end{pmatrix}
=\avec(t,\tmsoe,\xxvec,\vvvec) 
+ \eps \avec^1 (t,\tmsoe,\xxvec,\vvvec)
+ \frac 1\eps \bvec(t,\tmsoe,\xxvec,\vvvec),
\end{equation}
very close to the one studied by Fr\'enod \cite{frenod2004}
which originates from gyrokinetic plasma questions 
(see also 
Poincar\'e \cite{poincareMNMC},
Krylov and Bogoliubov \cite{KryBo},
Bogoliubov and Mitropolsky \cite{BogoMitro},
Sanders and Verhulst \cite{SanVer},
Schochet \cite{schochet:1994},
Fr\'enod and Sonnendr\"ucker \cite{frenod/sonnendrucker:1997,
frenod/sonnendrucker:1998, frenod/sonnendrucker:1999},
Fr\'enod and Watbled \cite{FW} and
Fr\'enod, Raviart and Sonnendr\"ucker \cite{FRS:1999}
for presentations of methods allowing to remove time
oscillations).
With very little changes we may apply the results it contains to deduce that
\begin{equation}\label{Rexp}
\begin{aligned}\ds 
\xxvec(t)
= \xxze(t,\tmsoe)
+ \eps\xxun (t,\tmsoe)
+ \dots, ~~~~~~
\vvvec(t)
= \vvvec^0(t,\tmsoe)
+ \eps\vvvec^1 (t,\tmsoe)
+ \dots,
\end{aligned}
\end{equation}
where oscillating functions $\xxze$, $\vvvec^0$, $\xxun$ and  $\vvvec^1$
are linked to non oscillating functions 
$\yyze$, $\uuvec^0$, $\yyun$ and  $\uuvec^1$
by 
\begin{align}
& \xxze(t,\theta)= \yyze(t),
\label{X0-Y0-1}
\\
& \vvvec^0(t,\theta)=\vmr(\tc,\theta,\yyze(t)) 
  + \uuvec^0(t),
\label{V0-U0-1}
\end{align}
and 
\begin{align}
&\xxun(t,\theta)= \yyun(t)
+\int_0^\theta\vmr(\tc,\sigma,\yyze(t))\, d\sigma,
\label{X1-Y1-1}
\\
&
\begin{aligned}
       \vvvec^1 &(t,\theta)=  
       \{\nabla\vmr(\tc,\theta,\yyze(t))\}
       \{\yyun(t)
         +\int_0^\theta\vmr(\tc,\sigma,\yyze(t))\, d\sigma\}
      \\ \ds &+\uuvec^1(t)
       + \pvmr(\tc,\theta,\yyze(t)) 
       - \pvmr(\tc,0,\yyze(t)) 
       + \int_0^\theta\Big(\vwn(\tc,\sigma,\yyze(t))- \bint\vwn(\tc,\varsigma,\yyze(t))\, d\varsigma\Big)\, d\sigma 
      \\ 
       &
        - \int_0^\theta \vmr(\tc,\sigma,\yyze(t))\, d\sigma.
 \end{aligned}
\label{U1-V1-1}
\end{align}
Then $\yyze$, $\uuvec^0$, $\yyun$ and  $\uuvec^1$ are the solution to
\begin{align}
&\frac{d\yyze}{dt} = \uuvec^0 ,
\label{eqY0}
\\
&\frac{d\uuvec^0}{dt} =\bint\vwn(\tc,\theta,\yyze)\, d\theta  -\uuvec^0, 
\label{eqU0}
\end{align}
\begin{align}
&\begin{aligned}
\frac{d\yyun}{dt}&=
       \ds \bint\{\nabla\vmr(\tc,\theta,\yyze)\}
       \{\int_0^\theta\vmr(\tc,\sigma,\yyze)\, d\sigma\} \,d\theta
      \\ \ds 
       &+\uuvec^1
       + \bint\pvmr(\tc,\theta,\yyze) \,d\theta
       -\pvmr(\tc,0,\yyze)
       + \bint\int_0^\theta\Big(\vwn(\tc,\sigma,\yyze)
       - \bint\vwn(\tc,\varsigma,\yyze)\, d\varsigma \Big) d\sigma\, d\theta
      \\
       &- \bint\int_0^\theta \vmr(\tc,\sigma,\yyze)\, d\sigma d\theta
       -\Big\{\bint\int_0^\theta
               \nabla\vmr(\tc,\sigma,\yze)\, d\sigma d\theta \Big\}\{\uuvec^0\}
       -\bint\int_0^\theta\fracp{\vmr}{\tc}(\tc,\sigma,\yze)\, d\sigma d\theta,
\end{aligned}
\label{eqY1}
\end{align}
\begin{align}
&\begin{aligned}
\frac{d\uuvec^1}{dt}&=
        \ds \big\{\bint\nabla\vwn(\tc,\theta,\yyze)d\theta\big\} \{\yyun\}
        +\bint\big\{\nabla\vwn(\tc,\theta,\yyze)\big\} 
          \big\{\int_0^\theta\vmr(\tc,\sigma,\yyze)\, d\sigma\big\}d\theta
      \\
       &-\Big\{\bint\{\nabla\vmr(\tc,\theta,\yyze)\}
       \{\int_0^\theta\vmr(\tc,\sigma,\yyze)\, d\sigma\} d\theta
       +\uuvec^1
       + \bint\pvmr(\tc,\theta,\yyze) \,d\theta
       - \pvmr(\tc,0,\yyze)
      \\ 
       &+ \bint(\int_0^\theta\vwn(\tc,\sigma,\yyze)\, d\sigma 
       - \theta\bint\vwn(\tc,\sigma,\yyze)\, d\sigma ) \,d\theta
       - \bint\int_0^\theta \vmr(\tc,\sigma,\yyze)\, d\sigma d\theta\Big\}
      \\\ds 
       &
       +\Big\{ \nabla\pvmr(\tc,0,\yyze)
     \\ \ds 
       & - \bint\nabla\Big(\int_0^\theta\vwn(\cdot,\sigma,\cdot)\, d\sigma  
       - \theta\bint\vwn(\cdot,\sigma ,\cdot)\, d\sigma \Big)(\tc,\yyze)\, d\theta
       + \bint\int_0^\theta 
           \vmr(\tc,\sigma,\yyze)\, d\sigma d\theta\Big\}\{\uuvec^0\}
     \\ \ds 
       & + \fracp{\pvmr}{t}(\tc,0,\yyze)
       - \bint \fracp{\ds\Big(\int_0^\theta\vwn(\cdot,\sigma,\cdot)\, d\sigma 
       - \theta\bint\vwn(\cdot,\sigma ,\cdot)\, d\sigma  \Big)}{t}(\tc,\yyze)\,d\theta
     \\ \ds & ~~~~~~~~~~~~~~~~~~~~~~~~~~~~~~~~~~~~~~~~~
              ~~~~~~~~~~~~~~~~~~~~~~~~~~~~~~~~~~~~~~~~~
       + \bint\int_0^\theta \fracp{\vmr}{t}(\tc,\sigma,\yyze)\, d\sigma d\theta,
\end{aligned}
\label{eqU1}
\end{align}
equipped with initial conditions
$\yyze(0)=\xvec$, $\uuvec^0(0)=\vvec$, $\yyun(0)=0$ and  $\uuvec^1(0)=0$.

~

As said previously, this result may be inferred from \cite{frenod2004}.
It may also be deduced by expanding fields and functions in 
(\ref{SM1ts1}) and (\ref{SM1ts2}) in a correct way. This is done formally in the Appendix.

~

Equations (\ref{X0-Y0-1}) and (\ref{V0-U0-1}) mean that the order 0 trajectory
does not oscillate and that the 0 order velocity is the tide velocity 
added with a non oscillating velocity $\uuvec^0$ which is generated by wind.
Wind acts on $\uuvec^0$ only through its averaged value. This is translated
in (\ref{eqU0}). By the way, notice that 
since the average value of $\vmr$ is $0$, equations (\ref{eqY0})-(\ref{eqU0})
only involve the average wind. 
As concerns order 1 terms, the situation is more complex. 
First since  $\int_0^\theta\vmr(\tc,\sigma,\yyze(t))\, d\sigma$ may be 
interpreted as the position of a sea water particle placed in $ \yyze(t)$
at the beginning of a tide cycle ($\theta=0$), equation (\ref{X1-Y1-1})
means that the order 1 position is this water particle position plus a non 
oscillating function $\yyun(t)$.
Regarding the terms (\ref{U1-V1-1}) contains, the first one describes the
way the space variation of the tide velocity acts, the second one is the
non oscillating part of the velocity.
The third  and fourth terms quantify the action of the sea velocity perturbation.
Concerning the next term, we need to remember that the action of the averaged value
of the wind velocity is taken into account in order 0 equation (\ref{eqU0}).
Then we notice that 
$\vwn(\tc,\sigma,\yyze(t))- \bintsml\vwn(\tc,\varsigma,\yyze(t))$
quantifies the wind action around its averaged value at each time of 
a tide cycle. Hence 
$\int_0^\theta\Big(\vwn(\tc,\sigma,\yyze(t))- 
\bintsml\vwn(\tc,\varsigma,\yyze(t)))\, d\varsigma\Big)\, d\sigma$
is the cumulated action of the wind around its averaged value. This quantity
acts at order 1.
The last term of (\ref{U1-V1-1}) can be interpreted as the previous one 
recalling that the mean value of the sea velocity is 0.
It is hard to give intuitive explanations for the evolution equations  
(\ref{eqY1})-(\ref{eqU1}). We only notice that they
involve mean value of non-linear interactions between fields which quantify
the mean joint action of sea and wind.
This non intuitive quantification is made possible thanks to the asymptotic analysis
presented in \cite{frenod2004} and in the Appendix.

The characteristic size of the right hand sides in (\ref{eqY0})-(\ref{eqU1})
are 1. Hence, to get the same precision
about $\eps^2$, with the same Euler scheme using (\ref{Rexp})-(\ref{eqU1})
to compute $(\xxvec, \vvvec)$, in place of using (\ref{SM1ts1})-(\ref{SM1ts2})
as evoked in the beginning of this section, the needed time step $\Delta t$
is about $\eps^2$. 
Comparing this with the value $\eps^3e^{-\frac{1}{\eps}}$ found in the 
beginning
of the section, despite the heavy form of (\ref{eqY0})-(\ref{eqU1}),
this is an appreciable gain even more so given that 
our probabilistic forecast method needs numbers of simulations
for numbers of wind time series.

~

Concerning the validity of the expansion (\ref{Rexp}),
if  $\vwn(\tc,\theta,\xvec)$ was a regular function, periodic 
with respect to $\theta$, 
asymptotic expansion (\ref{Rexp}) could be rigorously justified, 
applying \cite{frenod2004}, by the following inequality
\begin{gather}
\sup_{t\in [0,1]}\|\xxvec(t)-\xxze(t,\tmsoe)\|\leq c\eps, ~ 
\sup_{t\in [0,1]}\|\xxvec(t)-\xxze(t,\tmsoe)- \eps\xxun (t,\tmsoe)\|\leq c\eps^2,
\\
\sup_{t\in [0,1]}\|\vvvec(t)- \vvvec^0(t,\tmsoe)\|\leq c\eps, ~ 
\sup_{t\in [0,1]}\|\vvvec(t)-\vvvec^0(t,\tmsoe)-  \eps\vvvec^1 (t,\tmsoe) \|
\leq c\eps^2,
\end{gather}
which would be true for $\eps$ small enough and  
for a constant $c$ independent of $\eps$, where $\|\cdot \|$ 
stands for the Euclidean norm in $\rit^2$.

~

As we cannot consider the wind time series to have this form,
we will estimate the error using numerical experiments. This is done in the next section: several
realistic wind time series are simulated using a method described in 
Monbet {\it et al.} \cite{MonbetAlliotPrevosto} 
and the corresponding error is calculated for each of them.
For this, we have to decide
how the average values are to be computed. If the way is clear
for fields linked to $\vmr$ or $\pvmr$, for which we simply
take $\bintsml d\theta  = \int_0^1  d\theta$, the situation is more
obscure for wind linked fields.
For them, we will take as averaged value at time $t$ the mean
value on the interval centered in $t$ and with length $p$, i.e.
\begin{gather}
\bint \vwn(\tc,\theta,\xvec)\,  d\theta = 
\frac{1}{p}\int_{t-p/2}^{t+p/2}\vwn(s,\frac s\eps,\xvec)\, ds,
\end{gather}
for a parameter $p$ that has to be adjusted experimentally.
As concerns the integrals in $\theta$ of the wind linked
fields, they have to be replaced by integrals respecting 
$\vvvec_1(t,\toe)=\uuvec_1(t)$ when $\toe$ is an integer,
which is a constraint imposed by the Averaging Method. 
Hence, we replace 
\begin{gather}
\bigg(\int_0^\theta \vwn(\tc,\sigma,\xvec ) 
\,d\sigma\bigg)_{|\theta=\toe} \text{ ~ by ~ }
\int_{\eps[\toe]}^{t} \vwn(s,\frac s\eps,\xvec)\, ds,
\end{gather}
where $[\toe]$ stands for the integer part of $\toe$.
Those choices have also to be experimentally validated.

\section{Numerical validation}
We begin this section by introducing the sea velocity fields
and briefly presenting the stochastic method used to generate the wind time series.
Then, the expansion (\ref{Rexp}) is validated by numerical experiments.
More precisely, several wind time series are generated. Then, for each 
of them, we compute the real trajectory using (\ref{SM1ts1})-(\ref{SM1ts2}) and compare it to the trajectory obtained with expansion (\ref{Rexp}) and the following equations.

\subsection{Metocean fields}

Let us first describe the metocean fields that have been used in 
our numerical experiments. Regarding the rescaled ocean velocity induced by the tide 
wave and its perturbation, we have chosen the following simple 
parametric fields:
\begin{gather}
 \vmr(t,\theta,\xvec) = (2+\sin(6\pi t)) ~ x_1 
 \begin{pmatrix}
\sin(2\pi\theta)+\frac{1}{4} \sin(4\pi\theta) \\ \frac{1}{2} \sin(2\pi\theta)
 \end{pmatrix},
\\
  \pvmr(t,\theta,\xvec) = (2+\cos(6\pi t)) ~ x_2 
  \begin{pmatrix} 
  \sin(2\pi\theta) \\ \sin(2\pi\theta)
  \end{pmatrix},
\end{gather}
where $x_1$ and $x_2$ are the first and second components of $\xvec$.

The trajectories associated with the field
$\begin{pmatrix}
\sin(2\pi\theta)+\frac{1}{4} \sin(4\pi\theta) \\ \frac{1}{2} \sin(2\pi\theta)
\end{pmatrix}$
are non circular loops which remind those we can observe in Nihoul \cite{Nihoul} or
Salomon and Breton \cite{SaBre}. 
The velocity field $\vmr$ is this simple vector field
modulated by a time and position dependency in order to see the influence
of time derivative and gradient on the object trajectory.  $\pvmr$ is also a
simple field which gradient is orthogonal to the one of  $\vmr$.

In order to simulate realistic wind time series, we have 
assumed that they can be decomposed as the sum of two components, 
e.g $\vwn(t,\toe,\xvec)=\vwn_{\hspace{-2pt}Lt}(t,\xvec)+\vwn_{\hspace{-2pt}st}(\toe,\xvec)$ where
\begin{itemize}
	\item $\vwn_{\hspace{-2pt}Lt}(t,\xvec)$ represents the wind evolution at a 
synoptic-scale, e.g at the scale of the high- and low-pressure systems. 
The typical dimension of these systems ranges approximately from 1000km 
and 2500km and their duration is a couple of days to at most a couple of weeks. 
	\item $\vwn_{\hspace{-2pt}st}(\toe,\xvec)$ represents the small scale 
evolution of the 
wind (e.g mesoscale and microscale winds). This scale includes phenomena 
such as thunderstorms, squall lines, land and sea breezes...
\end{itemize}

Such a decomposition is discussed more precisely by Breckling 
(see Breckling \cite{Breckling} and \cite{AilliotThese}).
Different methods have been proposed in the literature to simulate 
realistic wind time-series at the synoptic scale 
(see Monbet, Ailliot and Prevosto 
\cite{MonbetAlliotPrevosto}
and references therein). In this study, we have first assumed that the wind 
$\vwn_{\hspace{-2pt}Lt}(t,x)$ is homogeneous in space; i.e. 
$\vwn_{\hspace{-2pt}Lt}(t,x)=\vwn_{\hspace{-2pt}Lt}(t)$ for all 
$\xvec$ and $t$, what seems realistic according the size of the domain
which is supposed to be about some hundred kilometers. 
Then, we have used a non-parametric resampling method to simulate the process 
$\vwn_{\hspace{-2pt}Lt}(t)$, namely the Local Grid Bootstrap 
algorithm proposed by Monbet, Ailliot and Prevosto
\cite{MonbetAlliotPrevosto}. 
This method has already been validated on various datasets, and it was found 
that it can successfully be used to simulate realistic wind time series. 
In the present study, it has been calibrated on a dataset which describes 
the wind condition during the summer at a point of coordinates 
$(46.25\,N,1.67\,W)$, located near the French 
Atlantic coast. It describes the synoptic wind conditions during the last 
20 years recorded every $\Delta t_1=6$~hours.

Then, in order to simulate the small scale variations, we have used an 
Autoregressive model. More precisely, for simplicity reasons, we 
have first assumed that this field is homogeneous in space. This assumption 
is unrealistic and could be refined later on. Then, we have assumed that, 
for $k \in \nit^*$, 
\begin{gather}
\vwn_{\hspace{-2pt}st}(k \Delta t_2)=
a\vwn_{\hspace{-2pt}st}((k-1) \Delta t_2)+E(k \Delta t_2)
\end{gather}
where $\{E(k \Delta t_2)\}_{k\in \nit^*}$ 
denotes a zero-mean Gaussian white noise with covariance 
matrix $\sigma^2 I_2$. In practice, we have simulated this small-scale 
component with a time-step $\Delta t_2 = \eps /100$. This consists in
giving the wind value every $4$~min approximately,
 and the parameters $a$ and $\sigma$ have been chosen such that the 
process $\vwn_{\hspace{-2pt}st}$ has a memory of a few hours and 
such that the standard 
deviation of its marginal distribution represents approximately 
$10\%$ of the one of the process $\vwn$.

Finally, the procedure described above makes it possible to 
simulate the processes 
$\vwn_{\hspace{-2pt}Lt}(t)$ and $\vwn_{\hspace{-2pt}st}(\toe)$ 
for $t\in k\Delta t_1$ and $\toe\in k\Delta t_2$, 
respectively. The values of these processes for other values of $t$ 
have been calculated by linear extrapolation. An example of simulated 
wind time-series is given in Figure \ref{fig:fig2}.

\subsection{Numerical results}

In order to validate the asymptotic expansion given in the previous section, 
we have computed the solutions of (\ref{SM1ts1})-(\ref{SM1ts2})
and of (\ref{Rexp})-(\ref{eqU1})  for $N=100$ 
artificial wind time series. Let us denote 
$\hat{\xxvec}_i$, $\hat{\vvvec}_i$,
$\hat{\xxvec}_i^0$, $\hat{\vvvec}_i^0$, 
$\hat{\xxvec}_i^1$, $\hat{\vvvec}_i^1$, 
for $i\in\{1,...,N\}$, the corresponding numerical 
approximations of $\xxvec$, $\vvvec$, $\xxvec^0$, $\vvvec^0$,
$\xxvec^1$ and $\vvvec^1$, respectively. We have tested different 
algorithms to solve these ODEs, and the best 
results have been obtained with an explicit Runge-Kutta (4,5) formula. 
In practice, we have used the Maltlab's function ode45 which  is based on an explicit Runge-Kutta (4,5) formula, the Dormand-Prince pair \cite{dor}. In the numerical 
results given hereafter, we have used $\eps=1/50$. This choice makes it 
possible to compute the solutions of the system (\ref{SM1ts1})-(\ref{SM1ts2}) 
with a good precision 
in a reasonable computational time, and permits also easier graphics 
representations. We have used $\xvec_i=(1,1)$ and $\vvec_i=(0,0)$ as initial values.


In Table \ref{tab:tab1}, the norms of error in 
object position and velocity for order 0 and 1 expansions are given. Let us discuss more precisely the results obtained with $\epsilon = 1/50$. It shows that the asymptotic 
expansions obtained with $p=\eps/10$ (corresponding to an interval
of $40$~min), 
$p=\eps/2$ (corresponding to an interval
of $6\textrm{h}15\textrm{min}$) and 
$p=\eps$ (corresponding to an interval
of $12\textrm{h}30\textrm{min}$) are about 
$5\eps^2$ worth, and close to each other. For comparison, 
we also compute the solution to  system (\ref{Rexp})-(\ref{eqU1}) 
when the wind is null, 
and we also found an error equal to $0.0022 \approx 5\eps^2$. 
For $p=4\eps$ (corresponding to an interval
of $2$~days), the error is significantly higher.  
As concerns the computational coast, it decreases as $p$ increases, 
so that a good compromise seems to use $\eps/2 \leq p \leq \eps$; 
i.e. the wind at a synoptic scale. Such a value of $p$ permits to compute 
the solutions of the system (\ref{Rexp})-(\ref{eqU1}) with a good precision 
in a computational 
time significantly lower than the one corresponding to the system 
(\ref{SM1ts1})-(\ref{SM1ts2}). For instance, if an Euler scheme is used, the exact system (\ref{SM1ts1})-(\ref{SM1ts2}) requires about 1000 more iterations than the approximate system (\ref{Rexp})-(\ref{eqU1}) to achieve the same accuracy. And, for the problem considered in this paper, an iteration of the exact system is 20 less expensive in terms of computational time. But this last remark is not general since the computational time highly depends on the nature of the tide and current fields $M$ and $N$ : here they are modeled by a quite simple analytical formula.

\begin{figure}
\centering
\makebox{\includegraphics{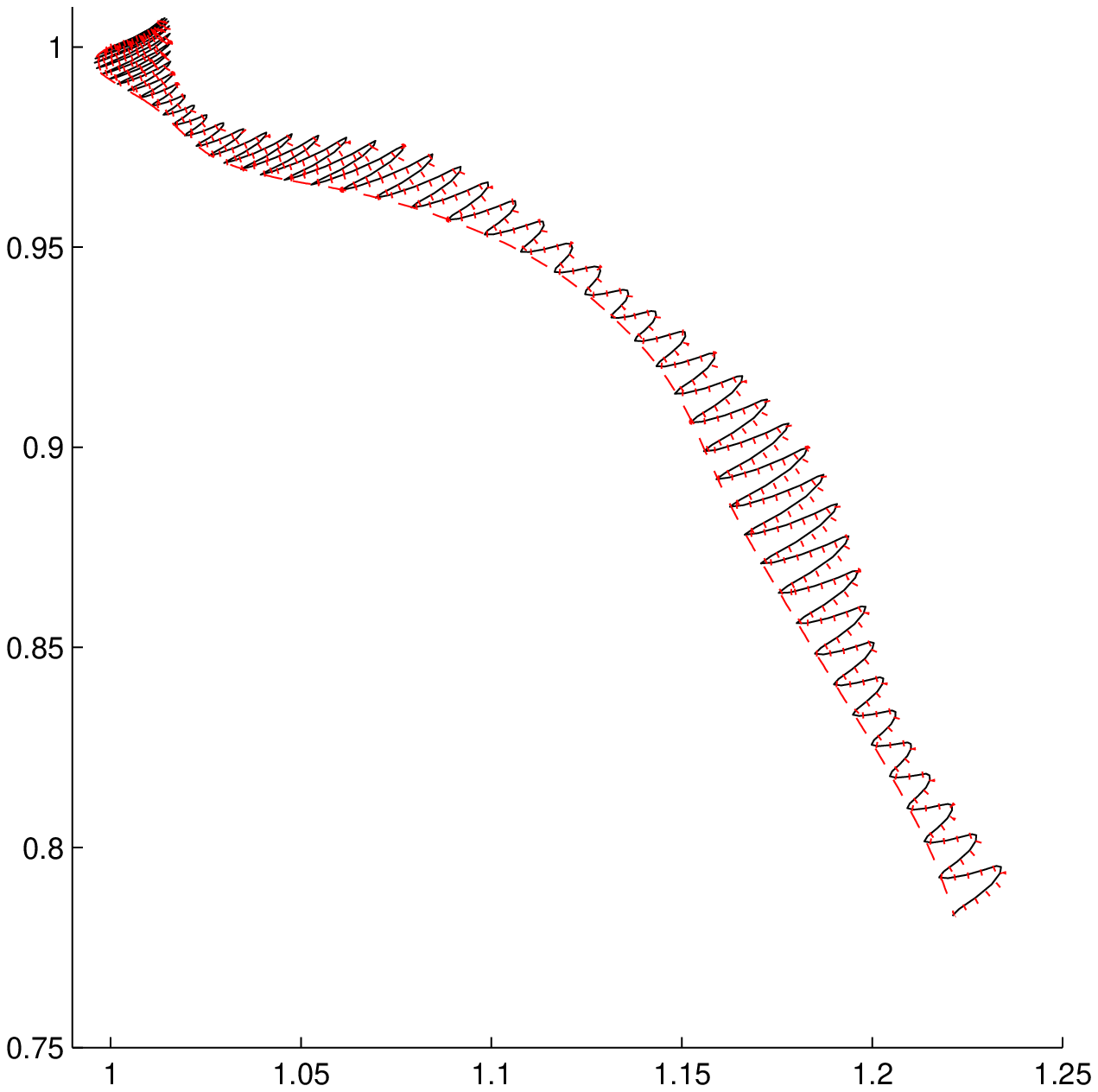}}
\caption{\label{fig:fig1} { Example of object's trajectory 
(two-dimensional phase plane plot). Solid line: $\hat{\xxvec}$, 
Dotted line: $\hat{\xxvec}^0+\eps \hat{\xxvec}^1$, 
Dashed line: $\hat{\yyvec}^0+\eps \hat{\yyvec}^1$}. $\eps = 1/50, p = \eps/2$.}
\end{figure}

\begin{figure}
\centering
\makebox{\includegraphics{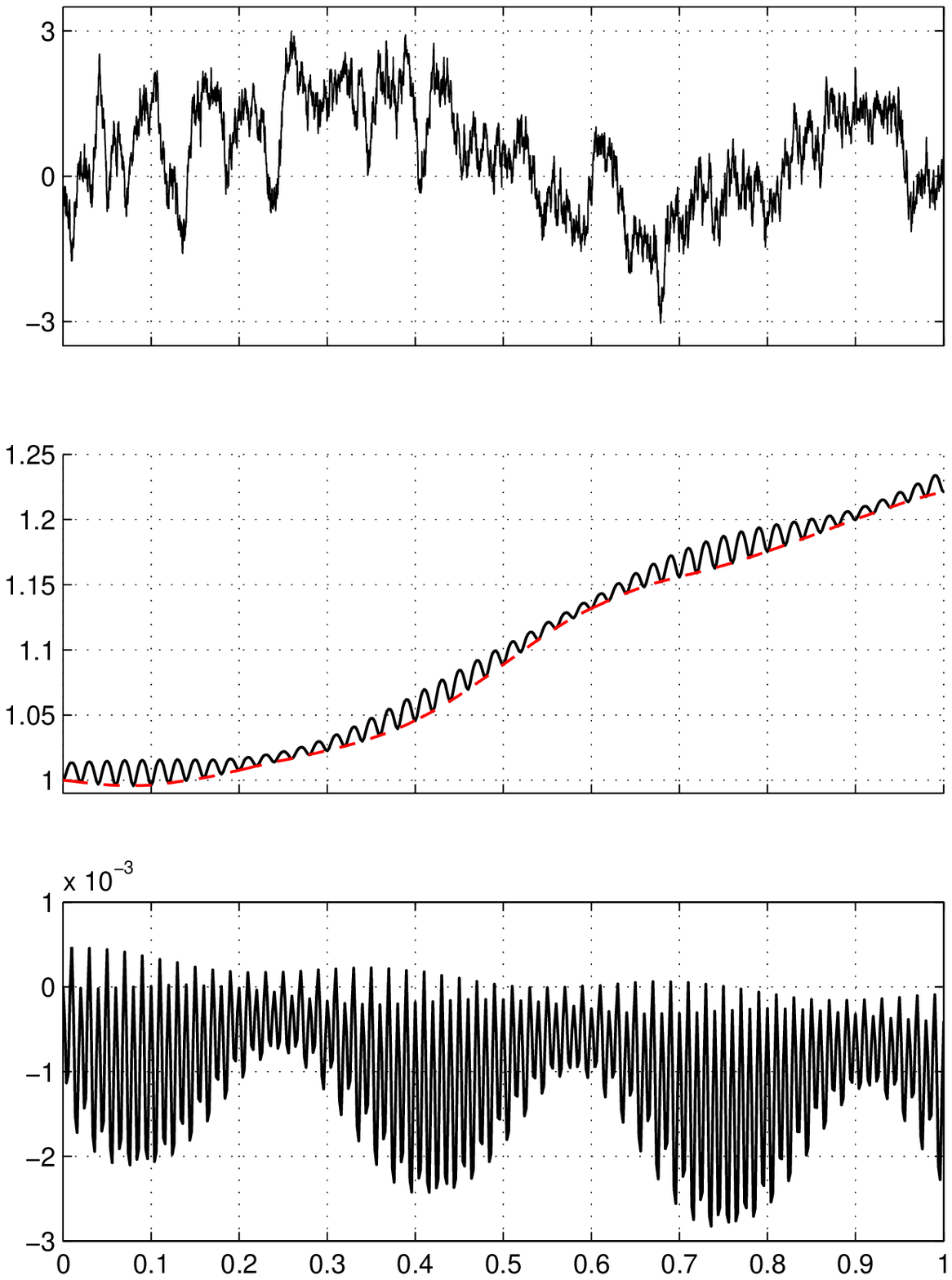}}
\caption{\label{fig:fig2} {Example of object's trajectory 
(time series plot)(Top: zonal wind first component, 
Second: first component of zonal object position 
(Solid line: $\hat{\xxvec}$, 
Dashed line: $\hat{\yyvec}^0+\eps \hat{\yyvec}^1$), 
Third: zonal error $\|\hat{\xxvec}-\hat{\xxvec}^0-\eps \hat{\xxvec}^1\|$). $\eps = 1/50, p = \eps/2$}}
\end{figure}

\begin{figure}
\centering
\makebox{\includegraphics{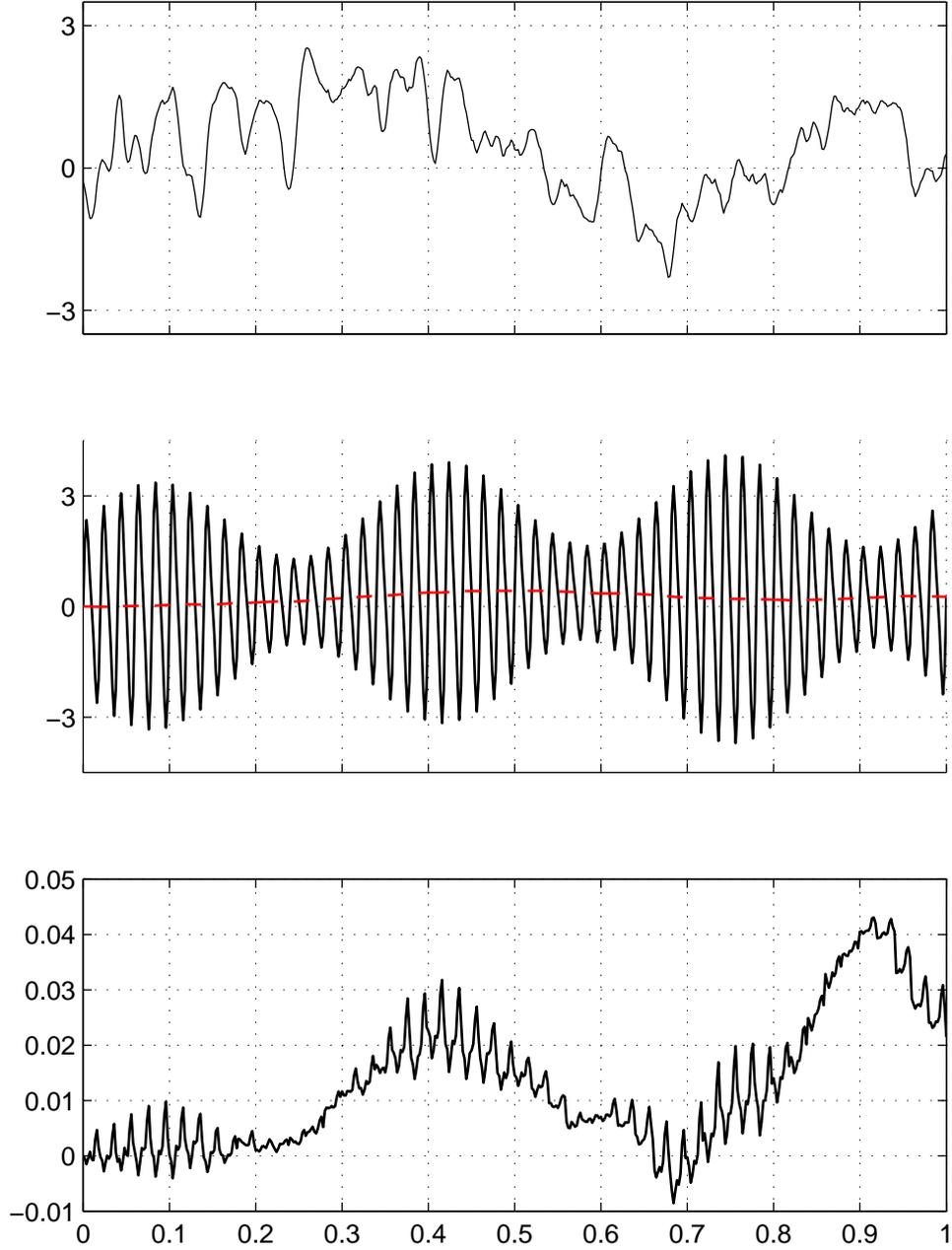}}
\caption{\label{fig:fig2-1} {Example of object's speed 
(time series plot)(Top: smoothed zonal wind component 
($p=\epsilon/2$), Second: first component of zonal object's speed 
(Solid line: $\hat{\vvvec}$, 
Dashed line: $\hat{\uuvec}^0+\eps\hat{\uuvec}^1$), 
Third: zonal error $\|\hat{\vvvec}-\hat{\vvvec}^0-\eps \hat{\vvvec}^1\|$). $\eps = 1/50, p = \eps/2$}}
\end{figure}

In Figure \ref{fig:fig1}, we have plotted the object trajectory associated
with the wind time series shown on the top of Figure \ref{fig:fig2}, 
using $p=\eps/2$. More precisely, the solid line represent $\hat{\xxvec}$
computed by directly solving (\ref{SM1ts1})-(\ref{SM1ts2}). The dashed 
line represents the average trajectory $\hat{\yyvec}^0+\eps \hat{\yyvec}^1$
obtained solving (\ref{eqY0}) and (\ref{eqY1}). We can see that this averaged
trajectory follows nicely the trend of $\hat{\xxvec}$. Then reconstructing
$\hat{\xxvec}^0$ and $ \hat{\xxvec}^1$ using (\ref{X0-Y0-1}) and (\ref{X1-Y1-1}),
we have represented the trajectory $\hat{\xxvec}^0+\eps \hat{\xxvec}^1$ 
in dotted line. The superimposition with $\hat{\xxvec}$ seems to be almost 
perfect. In order to go further in the result analysis, we turn to 
Figure \ref{fig:fig2}. On the top of this figure, we have shown the wind as
a function of time.  The second figure represents, in dashed line, the average trajectory
of the position of the first component and, in solid line, the first component
of the trajectory itself. We can see on this trajectory not only the long term trend
but also the two time periodic phenomena in game: tide oscillations 
(rapid oscillations) and tide coefficient amplitude (modulated amplitude).
Finally, the last plot exhibit $\eps^2-$order error 
(here $\eps =1/50=2. \, 10^{-2}$, then $10^{-3}=2.5 \, \eps^2$).
Moreover, we can see that this error function is a periodic function
with modulated amplitude. This spurs us thinking that we could improve the
accuracy of the reconstructed trajectory, if it is needed, considering the
next terms in the expansion of $\xxvec$ and $\vvvec$.
Figure \ref{fig:fig2-1} shows the first component of the average wind,
the average velocity (dashed line) and the velocity itself (solid line).
The third plot shows the error on the velocity. This figure permits
to visualize the action of the wind on the averaged velocity, which first
increases and then decreases.

~

In order to make our numerical validation more convincing, we also
made simulations for varying $\eps$. In Table \ref{tab:tab2}, the norms 
of the error on 
object position and velocity for order 0 and 1 expansions are 
given for $\eps = 1/100$ and $\eps = 1/25$ . The errors are proportional 
to $\eps$ for order 0
and to $\eps^2$ for order 1 as it was expected from the theory. 

\begin{table}
\centering
\fbox{%
\begin{tabular}{|l|cccc|}
  $\eps = \frac{1}{50}$& Speed (order 0) & Speed (order 1) & Position (order 0) & Position (order 1)\\
\hline

$p=\frac{\eps}{10}$    & 0.1115 [0.1115,0.1116]& 0.0460 [0.0454,0.0467] & 0.0174 [0.0173,0.0174] & 0.0023 [0.0021,0.0029]\\
$p=\frac{\eps}{2}$     & 0.1118 [0.1113,0.1134]& 0.0460 [0.0448,0.0473] & 0.0174 [0.0173,0.0174] & 0.0024 [0.0021,0.0032]\\
$p=\eps$  & 0.1125 [0.1108,0.1168]& 0.0460 [0.0441,0.0488] & 0.0173 [0.0171,0.0174] & 0.0026 [0.0020,0.0047]\\
$p=4\eps$  & 0.1192 [0.1102,0.1280]& 0.0542 [0.0438,0.0691] & 0.0174 [0.0166,0.0196] & 0.0065 [0.0024,0.0144]\\
$W\equiv0$ & 0.1115 & 0.0096& 0.0174&0.0022\\
\end{tabular}}

\caption{ \label{tab:tab1} 
Mean value, minimum and maximum values (mean [min, max]) of the errors $\sup_{t\in[0,1]}\|\hat{\vvvec}_i(t)-\hat{\vvvec}_i^0(t)\|$ (second column), 
$\sup_{t\in[0,1]}\|\hat{\vvvec}_i(t)-\hat{\vvvec}_i^0(t)-\eps \hat{\vvvec}_i^1(t)\|$ (third column),
$\sup_{t\in[0,1]}\|\hat{\xxvec}_i(t)-\hat{\xxvec}_i^0(t)\|$ (forth column) and
$\sup_{t\in[0,1]}\|\hat{\xxvec}_i(t)-\hat{\xxvec}_i^0(t)-\eps \hat{\xxvec}_i^1(t)\|$ (fifth column) for different values of $p$. 
The last line for $\eps = \frac{1}{50}$ gives the error for a zero wind field.}
\end{table}

\begin{table}

\centering
\fbox{%
\begin{tabular}{|l|cccc|}
\hline
 $\epsilon = \frac{1}{100}$ & Speed (order 0) & Speed (order 1) & Position (order 0) & Position (order 1)\\
\hline

$p=\frac{\eps}{10}$    & 0.0557 [0.0557,0.0557] & 0.0292 [0.0270,0.0358] & 0.0086 [0.0085,0.0086]& 0.0011 [ 0.0011,0.0014]\\
$p=\frac{\eps}{2}$     &0.0558 [0.0554,0.0563] & 0.0294 [0.0270,0.0362]& 0.0086 [0.0085,0.0086]& 0.0012 [0.0010,0.0015] \\
$p=\eps$  & 0.0560 [0.0552,0.0574] & 0.0296 [0.0273,0.0363] & 0.0085 [0.0084,0.0086]& 0.0013 [0.0010,0.0020]\\
$p=4\eps$  &0.0595 [0.0552,0.680] &0.0355 [0.0304,0.0516]& 0.0086 [0.0082,0.0098]& 0.0026 [0.0010,0.0072] \\

\hline
\hline
$\eps = \frac{1}{25}$ & Speed (order 0) & Speed (order 1) & Position (order 0) & Position (order 1)\\
\hline
$p=\frac{\eps}{10}$    &0.2254 [0.2249, 0.2257] &  0.1143 [0.1022 ,0.1418]&  0.0358  [0.0358,0.0359] & 0.0046 [0.0040,0.0057]\\
$p=\frac{\eps}{2}$     & 0.2254 [0.2226,0.2277] & 0.1142 [0.1017,0.1410] & 0.0358  [0.0357,0.0359]  & 0.0045 [0.0039,0.0058]\\
$p=\eps$  & 0.2257 [0.2201,0.2310]&  0.1142 [0.1011,0.1379]& 0.0358  [0.0356,0.0360] & 0.0044 [0.0037,0.0069]\\
\hline

\end{tabular}}
\caption{ \label{tab:tab2} 
Mean value, minimum and maximum values (mean [min, max]) of the errors $\sup_{t\in[0,1]}\|\hat{\vvvec}_i(t)-\hat{\vvvec}_i^0(t)\|$ (second column), 
$\sup_{t\in[0,1]}\|\hat{\vvvec}_i(t)-\hat{\vvvec}_i^0(t)-\eps \hat{\vvvec}_i^1(t)\|$ (third column),
$\sup_{t\in[0,1]}\|\hat{\xxvec}_i(t)-\hat{\xxvec}_i^0(t)\|$ (forth column) and
$\sup_{t\in[0,1]}\|\hat{\xxvec}_i(t)-\hat{\xxvec}_i^0(t)-\eps \hat{\xxvec}_i^1(t)\|$ (fifth column) for different values of $p$. }
\end{table}

\section{Application: long term forecast of an object's drift with Monte Carlo Method}

In the previous section, we have shown that the expansion 
of the trajectories of an object in near coastal ocean submitted to wind
makes it possible for us to compute quickly good approximation 
of this trajectory. It is then possible to make this computation for a 
large number of synthetic wind time series, and thus deduce the 
probability of a given scenario (this may be the probability of presence in a given area,
of collision or of running aground, for example). 

\begin{figure}
\centering
\makebox{\includegraphics{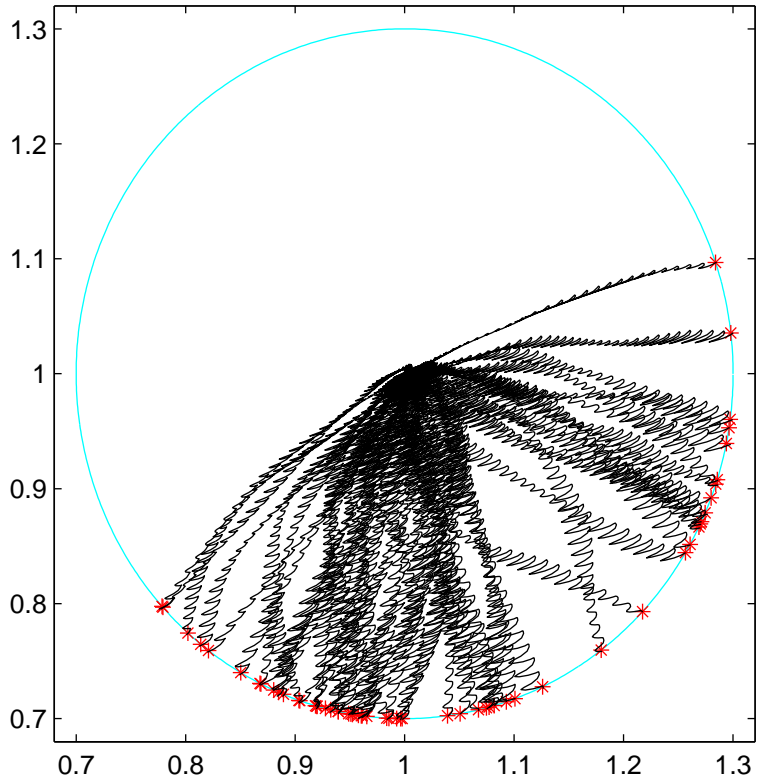}}
\caption{\label{fig:fig3} {Object's trajectory corresponding to 
100 synthetic wind time series and corresponding running aground points(*) }}
\end{figure}

\begin{figure}
\centering
\makebox{\includegraphics{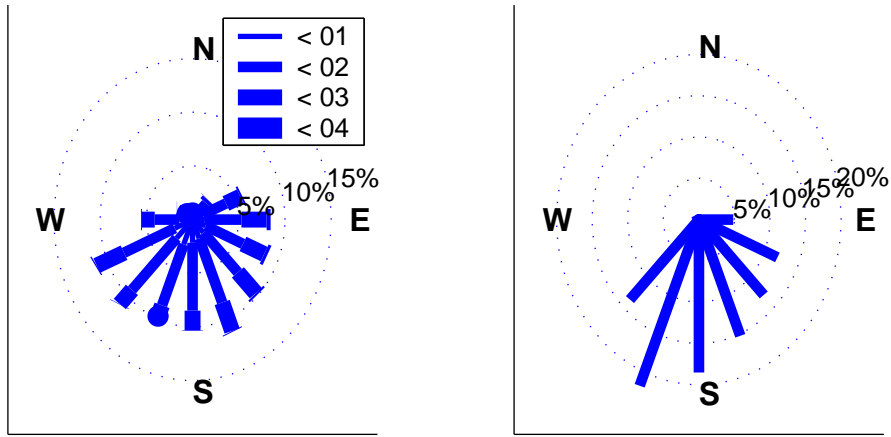}}
\caption{\label{fig:fig4} {Repartition of the wind (left) and of 
the running aground points (right) }}
\end{figure}

As an example, in this section we compute a running aground probability.
More precisely,  we have fixed an arbitrary coast 
line, e.g. the circle of center $(1,1)$ and radius $0.3$ and we have
focused on the running aground of the object on this coast. 
We have generated one thousand wind time series with the same stochastic model as in the
previous section and have also used the same sea velocity fields. We have computed
the thousand associated trajectory.
For this, we have used the method described in the previous sections,
based on asymptotic expansion (\ref{Rexp}).
The methodology is illustrated in Figure \ref{fig:fig3} 
where one hundred of the
thousand trajectories are drawn.
Proceeding this way, we have obtained enough trajectories 
to get reliable estimates
of the quantities of interest. First we can 
deduce an estimation of the probability of running aground 
by calculating the percentage of trajectories that reach the 
circle in the time interval $[0,1]$. 
We found $64.3\%$. Then we can estimate the 
probability that the object runs aground in a given area. 
The obtained results are shown in Figure \ref{fig:fig4}. On the left, 
we have plotted a wind rose which shows the joint repartition of 
direction and intensity. The length of each bar is proportional to the
proportion of wind having a direction around the direction the bar indicates.
The bars are then shared into four subbars which lengths are proportional 
to the proportion of wind having the corresponding magnitude.
We can see that the wind is generally blowing from the north.
The second histogram must be read in the following way. The length of each
bar is proportional to the proportion of trajectories that run aground
around the point of the coast which is in the bar direction.
We then see that the proportion of running aground trajectories 
has a maximum around the South-South-West of the coast line.

\section{Conclusion}
In this paper, we have explored a method to make probabilistic forecasts
of long term drift, in near coastal ocean, of object submitted to tide and
wind. The method is based on generating a large number of realistic wind
time series and, for each of them, to compute the associated 
object trajectory. 
The probabilities of trajectory linked events may then be estimated using 
this large number of trajectory realizations.
In order to achieve this goal, we have used a method proposed 
by Monbet, Ailliot and Prevosto  
\cite{MonbetAlliotPrevosto} to generate realistic wind time series
and the Averaging Method based on an asymptotic expansion of
trajectory of Fr\'enod \cite{frenod2004} to remove tide 
oscillations and compute quickly the object's trajectory.
This method has allowed us to compute the probability for an object to run
aground on an academic domain.

~

The quality of the results are good enough for us to contemplate going further
in exploring this new methodology.
Among the things to do in order to achieve the objective of computing
operationally events probabilities in real coastal areas we wish to broach the
following.

First, we will present a complete scale analysis of a coupled system
Shallow Water Equations - Newton Principle modeling the joint running 
of ocean and drifting object. The aim is to identify several regimes under
interest such as ``storm regime in coastal zone'' or 
``stillness above continental shelf'' and so on, and the corresponding 
equations.
This would give a way to compute the sea velocity fields $\vmr$ and $\pvmr$,
$\vmr$ being exclusively due to tide wave and $\pvmr$ to meteorological
factors. Notice also that $\pvmr$ could be divided into two parts:
one linked to pressure variations and another linked to wind.
Those aspects would then have to be incorporated in a software in 
order to access the fields $\vmr$ and $\pvmr$. 
Then we shall study the abilities of the method presented in the present 
paper to fit to pseudo-periodic fields $\vmr$ and $\pvmr$.
We are optimistic because the method already
works well with pseudo-periodic wind fields. The realism of the 
synthetic wind time series used as input of the numerical model 
could also be improved. In particular, the space-time model 
proposed by Ailliot, Monbet and Prevosto \cite{Ailliot2006} 
could be used to simulate non-homogeneous wind fields. 
Finally, the method could be implemented on real coastal areas.

~
\newline
{\bf Acknowledgements}

We would like to thank the referees for their pertinent comments and suggestions.
We would also like to thank Junko Murakami who has kindly accepted the
task of proofreading the manuscript and improving our poor english.
\appendix
\section{Appendix : Computations leading (\ref{X0-Y0-1}) - (\ref{eqU1})}
In this Appendix we formally lead the asymptotic expansion that gives equations
(\ref{X0-Y0-1}) - (\ref{eqU1}).

First, if we define $\yyvec(t)$ and $\uuvec(t)$ so that
\begin{align}
& \xxvec(t)= \yyvec(t),
\label{X-Y-1}
\\
& \vvvec(t)=\vmr(\tc,\tmsoe,\yyvec(t)) + \uuvec(t),
\label{V-U-1}
\end{align}
and inserting those in (\ref{SM1ts1}) and (\ref{SM1ts2}), we deduce 
\begin{align} \ds
\label{eqYYY}
&\frac{d\yyvec}{d\tc} =\vmr(\tc,\tmsoe,\yyvec) + \uuvec, 
\\
\label{eqUUU}
&\begin{aligned}
  \frac{d\uuvec}{d\tc} =\fracp{\pvmr}{\theta}(\tc,\tsep,\yyvec)
   + \eps \fracp{\pvmr}{\tc}(\tc,\tsep,\yyvec)
   + \eps\big( \nabla \pvmr (\tc,\tsep,\yyvec)\big)\big(\vmr(\tc,\tmsoe,\yyvec) + \uuvec \big)
   ~~~~~~~~~~~~&\\
   + \vwn(\tc,\tsep,\yyvec)-\vmr(\tc,\tmsoe,\yyvec) - \uuvec.&
 \end{aligned}
\end{align}

Then, we assume that $\yyvec$ and $\uuvec$ may be expanded in the following way:
\begin{align} \ds
\label{expYYYY}
&\yyvec(t) = \yyze(t) + \eps (\yyun(t)+\aavec^1(\tc,\tsep,\yyze(t)))
     + \eps (\yyvec^2(t)+\aavec^2(\tc,\tsep,\yyze(t),\yyun(t)))+ \cdots,
\\
\label{expUUUU}
&\uuvec(t) = \uuvec^0(t)+ \eps (\uuvec^1(t)+\bbvec^1(\tc,\tsep,\yyze(t)))
    + \eps (\uuvec^2(t)+\bbvec^2(\tc,\tsep,\yyze(t),\yyun(t)))+ \cdots,
\end{align}
for functions $\aavec^0(\tc,\theta,\yyze)$, 
$\aavec^1(\tc,\theta,\yyze,\yyun)$,  
$\bbvec^0(\tc,\theta,\yyze)$ and  
$\bbvec^1(\tc,\theta,\yyze,\yyun)$ 
being $1-$ periodic with respect to $\theta$ to be defined latter.
Using those asymptotic expansions in (\ref{X-Y-1}) and (\ref{V-U-1}) yields
(\ref{X0-Y0-1}) and (\ref{V0-U0-1}). Using again (\ref{expYYYY}) and (\ref{expUUUU})
in (\ref{eqYYY}) and (\ref{eqUUU}), expanding every functions using 
a Taylor expansion  yields
\begin{multline}
\label{eqexpYYY}
\frac{d\yyze}{dt}
+ \eps \frac{d\yyun}{dt}
+ \eps \fracp{\aavec^1}{t}(\yyze)
+ \fracp{\aavec^1}{\theta}(\yyze)
+ \eps\big(\nabla \aavec^1(\yyze)\big)\big(\frac{d\yyze}{dt}\big)
+ \eps \fracp{\aavec^2}{\theta} (\yyze,\yyun) 
+\cdots\, ,
\\
= 
\vmr(\yyze)
+\eps\big(\nabla \vmr(\yyze)\big)\big(\yyun+\aavec^1(\yyze)\big)
+\uuvec^0 + \eps (\uuvec^1+\bbvec^1(\yyze))
+\cdots
\end{multline}
\begin{multline}
\label{eqexpUUU}
\frac{d\uuvec^0}{dt}
+ \eps \frac{d\uuvec^1}{dt}
+ \eps \fracp{\bbvec^1}{t}(\yyze)
+ \fracp{\bbvec^1}{\theta}(\yyze)
+ \eps\big(\nabla \bbvec^1(\yyze)\big)\big(\frac{d\yyze}{dt}\big)
+ \eps \fracp{\bbvec^2}{\theta} (\yyze,\yyun) 
+\cdots
\\
=\fracp{\pvmr}{\theta}(\yyze)
+\eps\big(\fracp{\nabla \pvmr}{\theta}(\yyze)\big)\big(\yyun+\aavec^1(\yyze)\big)
+\eps\frac{d\pvmr}{dt}(\yyze)
+\eps\big(\nabla \pvmr(\yyze)\big)\big(\vmr(\yyze)+\uuvec^0\big)~~~~~~~~~
\\
+\vwn(\yyze)
+\eps\big(\nabla\vwn(\yyze) \big)\big(\yyun+\aavec^1(\yyze)\big)
-\vmr(\yyze)
-\eps\big(\nabla\vmr(\yyze) \big)\big(\yyun+\aavec^1(\yyze)\big)~~~~~~~~~
\\
-\uuvec^0
-\eps(\uuvec^1+\bbvec^1(\yyze))
+\cdots\, ,
\end{multline}
where $+\cdots$ contains every terms being of order greater than 1 in $\eps$.

Identifying the terms of order 0 in $\eps$ gives
\begin{align}
\label{eqY0000}
&\frac{d\yyze}{dt}
+\fracp{\aavec^1}{\theta}(\yyze)
=\vmr(\yyze)+\uuvec^0 ,
\\
\label{eqU0000}
&\frac{d\uuvec^0}{dt}
+\fracp{\bbvec^1}{\theta}(\yyze)
=\fracp{\pvmr}{\theta}(\yyze)
+\vwn(\yyze)
-\vmr(\yyze)
-\uuvec^0.
\end{align}
Integrating those equations with respect to $\theta$ from $0$ to $1$,
we obtain (\ref{eqY0}) and (\ref{eqU0}). Once this is done, 
equations (\ref{eqY0000}) and (\ref{eqU0000}) yield
\begin{gather}
\aavec^1(\tc,\theta,\yyze)= \int_0^\theta \vmr (\tc,\sigma,\yyze(t))\,d\sigma,
\\
\begin{aligned}
\bbvec^1(\tc,\theta,\yyze)&= \uuvec^1(t)
       + \pvmr(\tc,\theta,\yyze(t)) 
       - \pvmr(\tc,0,\yyze(t)) 
\\
       &+\int_0^\theta\Big(\vwn(\tc,\sigma,\yyze(t))- \bint\vwn(\tc,\varsigma,\yyze(t))\, d\varsigma\Big)\, d\sigma 
        - \int_0^\theta \vmr(\tc,\sigma,\yyze(t))\, d\sigma,
\end{aligned}
\end{gather}
and then (\ref{X1-Y1-1}) and (\ref{U1-V1-1}).

Looking now at the terms of order $1$ in $\eps$ in (\ref{eqexpYYY}) and 
(\ref{eqexpUUU}), after replacing $\aavec^1$, $\bbvec^1$, $d\yyze/dt$  
and $d\uuvec^0/dt$ by their
expressions, results in two equations. Those equations, when integrated in $\theta$,
lead in a heavy but straightforward way (\ref{eqY1}) and (\ref{eqU1}).

\bibliographystyle{plain}
\bibliography{biblio}

\end{document}